\documentclass[12pt]{amsart}
\usepackage{amssymb}
\usepackage{euscript}

\let\cal\mathcal

\input xypic
\xyoption {all}

\newtheorem{theorem}{Theorem}

\newtheorem{proposition}[theorem]{Proposition}

\newtheorem{definition}[theorem]{Definition}

\def\N{\mathbb N}

\def\M{\cal{M}}

\let\<\langle

\numberwithin{equation}{section}

\let\\\cr

\define\commadots{,.\,.\,.\,,}
\begin{document}

\title[Compact range property]{Compact Range Property and
Operators on $\boldsymbol C^{\boldsymbol*}$-algebras}
 
\author{Narcisse Randrianantoanina}
\thanks{Supported in part by NSF Grant DMS-9703789}
\address{Department of Mathematics and Statistics, Miami University, Oxford,
Ohio 45056}
\email{randrin@muohio.edu}
\subjclass{46L50,47D15}
\keywords{$C^*$-algebras, Vector measures}

\begin{abstract}
We prove that a Banach space $E$ has the 
compact range property (CRP) if
and only if for any given $C^*$-algebra $\cal A$, 
every  absolutely summing operator from
$\cal A$ into $E$ is compact. Related results for
$p$-summing operators ($ 0<p<1$) are also discussed as well
as operators on non-commutative $L^1$-spaces and 
$C^*$-summing operators. 
\end{abstract}

\maketitle

\section{INTRODUCTION}

A Banach space $E$ is said to have the compact range
 property (CRP)
if  every $E$-valued countably additive measure 
of bounded variation has compact range. 
 It is well known that every Banach space with
 the Radon-Nikodym property (RNP) 
has the (CRP) and for dual Banach spaces,
the (CRP) were completely characterized as
those whose predual do not contain  any copies $\ell^1$.  
For more in depth discussions on Banach
spaces with the (CRP), we refer to 
\cite{T4}.

The following characterization can be found in \cite{T4}:
  A Banach space $E$ has the (CRP) if and only if every
$1$-summing operator from $C[0,1]$ into $E$ is compact.  Since $C[0,1]$ is
a (commutative) $C^*$-algebra, it is a natural question whether $C[0,1]$
can be replaced by any $C^*$-algebras.  
Let us recall that in \cite{Ran7}, it 
was shown that if $X$ is a Banach space that
does not contain any copies of $\ell^1$ then
any $1$-summing operators from any given
$C^*$-algebra into $X^*$ is compact; hinting that, 
as  in commutative case, the (CRP) is the 
right condition to provide compactness.
The present note is an improvement of 
\cite{Ran7}. Our main result  confirms that, if $\cal
A$ is a $C^*$-algebra and $E$ 
is a Banach space that has the (CRP) then every $1$-summing from
$\cal A$ into $E$ is compact. Our proof relies on factorizations of
summing operators used in \cite{Ran7} and properties of
integral operators. 

There is another well kown characterization of 
spaces with the (CRP) in terms of operators defined on 
$L^1[0,1]$: a Banach space $E$ has the $(CRP)$ if and only
every operator $T$ from $L^1[0,1]$ into $E$ is Dunford Pettis
(completely continuous) thus the (CRP) is also referred to 
as the {\it complete continuity property} (CCP). 
Unlike the $1$-summing operators on $C^*$-algebras, 
operators defined
on non-commutative $L^1$-spaces do not behave the same way
as those defined on $L^1[0,1]$ do. In the last section of this
note, we will  discuss these operators along with 
$C^*$-summing operators studied by Pisier in \cite{PIS3}.

   Our
terminology and notation are standard as may be found in 
\cite{D1} and
\cite {DU} for Banach spaces, \cite {KR} and \cite {TAK} for $C^*$-algebras
and operator  algebras.

\section{PRELIMINARIES}

In this section, we recall some definitions.

\begin{definition}
Let $X$ and $Y$ be Banach spaces
and $0<p<\infty$.  An operator $T : X \to Y$ is said to be
$p$-summing if there is a constant $C$ such that for any finite sequence
$(x_1, x_2\commadots x_n)$ of $X$, one has
$$
\left(
\sum^n_{i=1} \Vert Tx_i\Vert^p\right)^{\frac1p}
 \leq C \sup\left\{\left(\sum^{n}_{i=1} \vert
\langle x_i, x^*\rangle \vert^p \right)^{\frac1p};
\  x^* \in X^*, \Vert x^*\Vert \leq 1 \right\} .
$$
\end{definition}

The smallest constant $C$ for which the above inequality holds is denoted
by $\pi_p(T)$ and is called the $p$-summing norm of $T$.

\begin{definition}
We say that an operator $T : X \to Y$ is an integral operator if it admits
a factorization:
$$
\xymatrix
{X \ar [r] ^{i\circ T} \ar[d]^\alpha &Y^{**}\\
L^\infty(\mu) \ar[r]^J &L^1(\mu) \ar[u]_\beta}
$$
where $i$ is the natural inclusion from $Y$ into $Y^{**},\ \mu$ is a
probability measure on a compact space $K$, $J$ is the natural inclusion
and $\alpha$ and $\beta$ are bounded linear operators.
\end{definition}

\noindent
We define the integral norm $i(T) : = \inf \left\{\Vert \alpha \Vert
\cdot\Vert \beta \Vert\right\}$ where the infinum is taken over all such
factorizations. 
 
 Similarily,  we shall say
that $T$ is {\it strictly integral} if $T$ is integral and on the factorization
above $\beta$ takes its values in $Y$.  

It is well known that integral operators are $1$-summing but
the converse is not true.

If $X = C(K)$ where $K$ is a compact Hausdorff space then
it is well known that 
every $1$-summing operator from $X$ into $Y$ is integral.

For more details on the different
properties of the classes of operators involved, we refer to \cite {DJT}.

The following simple fact will be needed in the sequel.
\begin{proposition}
  Let $T : X \to Y$ be a strictly integral operator.  If $Y$
has the (CRP) then $T$ is compact.
\end{proposition}
\begin{proof}
The operator $T$ has a factorization 
$T=\beta\circ J \circ \alpha$ where $\alpha: X \to L^\infty(\mu)$,
$J: L^\infty(\mu) \to L^1(\mu)$ and $\beta: L^1(\mu) \to Y$ are
as in the above definition.
Note that $J$ is $1$-summing so 
$\beta\circ J: L^\infty(\mu) \to Y$
is $1$-summing and since $L^\infty(\mu)$ is a $C(K)$-space and
$Y$ has the (CRP), $\beta\circ J$ (and hence $T$) is compact.
\end{proof}

\noindent
We recall that a von Neumann algebra $\cal M$ is said to be $\sigma$-finite
if the identity is countably decomposable equivalently if there exist a
faithful state $\varphi \in \cal M_*$.  As is customary, for every
functional $\varphi \in \cal M_*$ and $x \in \cal M$, $x \varphi$ (resp.
$\varphi x)$ denotes the normal functional $y \to \varphi (yx)$ (resp. $y
\to \varphi (xy))$.

\section{MAIN RESULT}

\begin{theorem}   For a Banach space $E$, the following are equivalent:
\begin{itemize}
\item [(1)] $E$ has the CRP;
\item [(2)]   Every $1$-summing operator $T: C[0,1] \to E$ is compact;
\item [(3)] For any given $C^*$-algebra $\cal A$, 
every $1$-summing operator $T : \cal{A}
\to E$ is compact.
\end{itemize}
\end{theorem}

\noindent
The equivalence $(1) \Leftrightarrow (2)$  is  well known, we
refer to \cite {DU}, \cite {T4} for more details.  Clearly (3)
$\Rightarrow$ (2)  so what we need to show is (1) 
$\Rightarrow$ (3).  For
this, it is enough to consider the following particular case 
(see \cite{Ran7} for this reduction).

\begin{proposition}  Let $E$ be a Banach space with the (CRP) 
and $\M$ be
a $\sigma$-finite von Neumann algebra.  If $T : \M \to E$ is $1$-summing and
is weak$^*$ to weakly continuous then $T$ is compact.
\end{proposition}

\begin{proof}  The proof is a refinement of the argument used in
Proposition~3.2 of \cite{Ran7}.  We will include most  of the details
for completeness. Without loss of generality, 
we can and do 
assume that $E$ is separable.

Let $\delta > 0$.  From Lemma~2.3 of \cite{Ran7},

$$
\Vert Tx\Vert \leq 2 (1 + \delta)
\pi_1(T)\Vert x f + fx \Vert_{\M_*} \quad \text{for every}
\  x \in  \M,
$$
\noindent
where $f$ is a faithful normal state in $\cal M_*$.  If $L^2(f)$ is
completion of the prehilbertian space $(\cal M, \langle \cdot, \cdot
\rangle)$ where 
$\displaystyle{\langle x,y\rangle = f(\frac {xy^* + y^*x} 2)}$
 then we
have the following factorization:

$$
\xymatrix
{
\cal M \ar[d]_J  \ar[r]^T & E\\
L^2(f) \ar[r]^\theta &L^2(f)^* \ar [r]^{J^*}&\M_* \ar [ul]_L
}
$$

\noindent
where $J$ is the inclusion map,
$\theta (Jx) = \langle \cdot, J(x^*)\rangle$ for every $x \in \cal M$ and
$L (\frac{xf+fx}{2}) = Tx$.
We recall that $L$ is a well defined bounded linear map since
$\{xf +fx; \ x \in \M \}$ is dense in $\M_*$ and
$\Vert L(xf +fx) \Vert \leq 4(1 +\delta) \pi_1(T) \Vert xf +fx \Vert_{\M_*}$.  
Let $ S := J^* \circ \theta \circ J$.

\medskip
\noindent
{\it Claim:  $J \circ L^* : E^* \to L^2(f)$ is compact.}
\medskip

  For this, let us
consider $L^*: E^* \longrightarrow \cal M$.
Since $E$ is separable, it is isometric to a subspace of
$C[0,1]$. Let 
$I_E$ be the isometric embedding of $E$ in $C[0,1]$ and $i$ 
be the natural
 inclusion of $C[0,1]$ into $C[0,1]^{**}$. 

 Define the following map
 $\widetilde T$ from $\cal M$ into $C[0,1]$ by setting
  $\widetilde T = \overline{I_E \circ T(x^*)}$ 
for every $x \in \M$.  (Here,
  $\overline f$ is the map $t \to \overline{f(t)}$ for 
$f \in C[0,1]$ with $\overline{f(t)}$ being the conjugate of the 
complex number $f(t)$).

Clearly, $\widetilde T$ is linear and bounded and 
it can be shown  that $\widetilde
T$ is $1$-summing and is weak$^*$ to weakly continuous.  In fact,
if $(x_1,x_2\commadots x_n)$ is a finite sequence in $\cal M$ 
then
\begin{equation*}
\begin{split}
\sum^n_{i=1} \Vert \widetilde T x_i \Vert 
&= \sum ^{n}_{i=1} \Vert
\overline{I_E \circ T(x_i^*)} \Vert\\
&= \sum^{n}_{i=1} \Vert I_E \circ T(x_i^*) \Vert\\
&= \sum^n_{i=1} \Vert T (x_i^*)\Vert\\
&\leq \pi_1(T)\sup\left\{\sum_{i=1}^n\vert\langle x_1^*, 
\varphi \rangle \vert,\
\varphi \in \cal M^*, \  \Vert \varphi\Vert \leq 1\right\}\\
&\leq \pi_1(T) \sup \left\{\sum_{i=1}^n \vert 
\langle x_i, \varphi^*\rangle\vert,\
\varphi \in \cal M^*, \  \Vert \varphi \Vert \leq 1
\right\}
\end{split}
\end{equation*}
so
$\widetilde T$ is $1$-summing with $\pi_1(\widetilde T) \leq \pi_1(T)$.
Moreover if $(x_\alpha)_\alpha$ is a net that converges to zero weak$^*$ in
$\M$ so does the net $(x^*_\alpha)_\alpha$ and since $T$ is weak$^*$ to
weakly continuous, $(T(x^*_\alpha))_\alpha$ converges to zero weakly in $E$
and hence  $(\widetilde T(x_\alpha))_\alpha$ is weakly null which shows
that $\widetilde T$ is weak$^*$ to weakly continuous.  

To complete the
proof, consider
$\xymatrix{
E^* \ar[r]^{L^*}&\cal M \ar [r]^-{i \circ \widetilde T} 
&C[0,1]^{**}.}
$

\noindent
Since $C[0,1]^{**}$ has the Hahn-Banach extension property
 and $i \circ
\widetilde T$ is $1$-summing, $i \circ \widetilde T$ 
is  an integral operator.
Let  $K : C[0,1]^* \longrightarrow \cal M_*$ such that 
$K^* = i \circ
\widetilde T$ (such operator exists since $i \circ \widetilde T$ is
 weak$^*$ to weakly continuous); $K$ is integral 
(\cite{DJT}) and since $\cal M_*$ is a complemented
subspace of its bidual $\M^*$ (see for instance \cite{TAK}), 
$K$ is strictly integral and therefore
$L\circ K : C[0,1]^* \to E$ is strictly integral and
by Proposition~3, $L\circ K$ (and hence
$(L\circ K)^* = i\circ\widetilde{T}\circ L^*$) 
is compact.

 Let  $(U_n)$ be a bounded sequence in $E^*$. There
 exists a subsequence $(U_{n_k})$ so that 
$(i \circ\widetilde{T}\circ  L^*(U_{n_k}))_k$ is
 norm convergent in $C[0,1]^{**}$.  
Since 
  $i$  and $I_E$ are  isometries, we get
 that $(T \circ L^*(U_{n_k}))_k$ is norm convergent so
 
$ \lim_{k,m} \Vert T(L^*(U_{n_k})^*) - T(L^*(U_{n_m})^*)\Vert = 0$.
 As in \cite{Ran7}, we get
$$
\lim_{k,m} \left\langle T(L^*(U_{n_k})^*)
- T(L^*(U_{n_m})^*), U_{n_k} - U_{n_m} 
\right\rangle = \lim_{k,m}\Vert
J \circ L^* (U_{n_k} - U_{n_m}) \Vert ^2_{L^2(f)}=0
$$
which proves that $\left(J \circ L^* (U_{n_k})\right)_k$ is norm-convergent in
$L^2(f)$.  The proof is complete.
\end{proof}


\begin{theorem}
Let $\cal{A}$ be a $C^*$-algebra, $E$ be a Banach space and
$0<p<1$. Every $p$-summing operator from $\cal A$ into $E$ is
compact.
\end{theorem}

\begin{proof}
Let $T: \cal{A} \to E$ be an operator with $\pi_p(T) <\infty$.
One can choose, by the Pietsch Factorization Theorem, a 
probability space $(\Omega, \Sigma, \mu)$ such that
$$
\xymatrix{
\cal{A} \ar[r]^T \ar[d]^J     &E  \\
 S  \ar[d] \ar[r]^{i_p} &S_p \ar[d] \ar[u]^{\widetilde T} \\
L^\infty(\mu) \ar[r]^{j_p} &L^p(\mu)
}$$
where $S$ is a subspace of $L^\infty(\mu)$, $S_p$ is the 
closure of $S$ in $L^p(\mu)$ and $i_p$ is the restriction of
the natural inclusion $j_p$.

Denote by $S_1$ the closure of $S$ in $L^1(\mu)$, by
$i_1$  the restriction of the natural inclusion and 
$i_{1, p}$  the natural inclusion of $S_1$ into $S_p$.

\medskip

\noindent
{\it Claim:  $\widetilde{T}\circ i_{1,p}: S_1 \to E$ is 
weakly compact.}

\medskip

To see this, let $(f_n)_n$ be a bounded sequence in
$S_1 \subset L^1(\mu)$. By Koml\`os's Theorem, there exists
a subsequence $(f_{n_k})_k$  and a function 
$f \in L^1(\mu)$ such that 
$\lim_{m \to \infty} \frac{1}{m} \sum_{k=1}^m  
f_{n_k}(\omega)= f(\omega)$ for a.e. $\omega \in \Omega$.
Since $0<p<1$, 
$$\lim_{m \to \infty} \Vert\frac{1}{m} \sum_{k=1}^m  
f_{n_k} - f \Vert_p =0.
$$ 
This shows that $f \in S_p$ and 
$ \left(\widetilde{T} \circ i_{1,p}\left(\frac{1}{m} \sum_{k=1}^m  
f_{n_k} \right) \right)_m$ converges to
$\widetilde{T}(f)$ in $E$ and the claim follows.

Using the factorization of weakly compact operator 
\cite{DFJP}, $i_{1,p}\circ\widetilde{T}$ factors through
a reflexive  space and since $i_1\circ J$ is $1$-summing, 
the theorem
follows from Theorem~4.
\end{proof}

\medskip

\section{Concluding remarks}
 
Let us recall some definitions
\begin{definition}
Let $X$  and $Y$ be  Banach spaces.  An operator
$T : X \to Y$ is called Dunford-Pettis if $T$ sends weakly compact
sets into norm compact sets.
\end{definition}

The following class of operators was introduced by 
Pisier in \cite{PIS3}
as extension of the $p$-summing operators in the setting of $C^*$-algebras.

\begin{definition}
Let $\cal{A}$ be a $C^*$-algebra and $E$ be a Banach space, $0<p<\infty$.
An operator $T:\cal{A}\to E$ is said to be $p$-$C^*$-summing if there exists a
constant $C$ such that for any finite sequence $(A_1,\ldots,A_n)$ of
Hermitian elements of $\cal{A}$,  one has
$$\left( \sum_{i=1}^n \|T(A_i)\|^p\right)^{\frac1p}
\le C\Big\| \left( \sum_{i=1}^n |A_i|^p\right)^{1/p}\Big\|_{\cal{A}}.$$
\end{definition}

Let $\M$ be a finite von Neumann algebra with a faithful
tracial state $\tau$ and let $J$ be the canonical
inclusion map from $\M$ into $L^1(\M,\tau)$. As in the commutative
case, we have the following :

\begin{proposition}
Let $ E$ be a Banach space and  $T: L^1(\M,\tau) \to E$ a 
bounded linear map. Then the following are equivalent:
\begin{itemize}
\item[(i)] $T$ is Dunford-Pettis;
\item[(ii)] $T\circ J$ is compact.
\end{itemize}
\end{proposition}

\begin{proof}
$(i) \implies (ii)$ is trivial. For the converse, let
$(a_n)_n$ be a weakly null sequence  in  the unit
ball of $L^1(\M,\tau)$.
It is clear that $(a_{n}^{*})_n$ is also weakly null so 
without  loss of generality, we can assume that 
$(a_n)_n$ is a sequence of self-adjoint operators.
For each $n\geq 1$, set 
$a_n= \int_{-\infty}^\infty t \ de_{t}^{(n)}$ the spectral
decomposition of $a_n$ and for every $N\geq 1$, let
$$p_{n,N}= \int_{-N}^N 1 \ de_{t}^{(n)}.$$
It is clear that for every $n\geq 1$ and $N \geq 1$,
$$\tau({\bf 1} - p_{n,N})=
\tau(\int_{\{|t| > N\}} 1 \ de_{t}^{(n)}) \leq \frac{1}{N}
\tau(|a_n|).$$
By the Akeman's characterization of relatively 
weakly compact subset in $L^1(\M,\tau)$
(see for instance \cite{TAK} Theorem~5.4 p.149), we conclude that
for any given $\epsilon >0$, there is $N_0 \geq 1$ such that
for every $n \geq 1$,
$\Vert a_n({\bf 1} - p_{n, N_0}) \Vert \leq \epsilon$. Moreover
$(a_n p_{n,N_0})_n$ is a bounded sequence in $\M$ and since
$T\circ J$ is compact, there is a compact subset
$K_\epsilon$ of $E$ such that $\{ T(a_n);\ n \in \N \} \subset
K_\epsilon + \epsilon B_E$. The proof is complete.

\end{proof}

Fix a type $II_1$ von Neumann algebra $\M$ such that $\M$
contains a complemented copy of a Hilbert space $H$. 
The space $H$ is reflexive (and  therefore has (CRP) ) but the
projection map $P$ from $L^1(\M,\tau)$ onto $H$ can not be 
Dunford-Pettis. 

A very well known property of $p$-summing operators is that 
they are Dunford-Pettis. This in not the case for
$C^*$-summing operators in general.
By Proposition~9, $P\circ J$ is not compact.
We remark that the argument used in \cite{Ran7} requires
 only that the operator is $C^*$-summing and Dunford-Pettis
hence since $J$ is clearly $C^*$-summing and $P\circ J$ is
not compact, $P\circ J$ should not be Dunford-Pettis. 

\bigskip

\noindent
{\bf Acknowledgments.} The author wishes to express his thanks to
Patrick Dowling  for helpful discussions regarding this note.

\providecommand{\bysame}{\leavevmode\hbox to3em{\hrulefill}\thinspace}

\end{document}